\documentstyle[12pt,amssymb,amscd]{amsart} 
\author{Dmitry E. Tamarkin}

\input amssym.def
\def \pt{\Pi T}
\def \g {{\frak g}}

\def \V{V[1]}
\def \proof {\begin{pf}}
\def \endproof {\end{pf}}
\def\gr{{\rm gr\;}}
\def \isom{\cong}
\def \pmin {\pm}
\def\mplus{\mp}
\def \G{{\cal G}}
\def \O{{\cal O}}

\def \F{{\cal F}}

\def \ass {{\cal A}ssoc}
\def \com {{\cal C}omm}
\def \lie {{\cal L}ie}
\def \e {e_2}
\def \bi{{\cal B}_{\infty}}
\def \b{{\cal B}}
\def \hass {{\cal H}oassoc}
\def \hcom {{\cal H}ocomm}
\def \hlie{{\cal  H}olie}
\def \he{{\cal HE}_2}
\def \a{{\cal A}}
\def \k{{\cal K}} 
\def \x {{\cal X}}
\def \p{P(\a)}
\def \t{{\cal T}^1_n} 
\def \cx{{\cal C}} 
\def\bbf{\Bbb}

\def \as{As}
\def \As{As}

\def\eqnref#1{(\ref {#1})}

\def \f{{\cal F}reeLie}

\def \co{{\rm Con}\;}
\def \cc {\co}
\def\v{{\cal V}}
\def \z {{\cal Z}}
\def \w {{\cal W}}
\newtheorem{Theorem}{THEOREM}[section]
\newtheorem{theorem}{THEOREM}[section]
 \newtheorem{Definition}[Theorem]{DEFINITION}
 \newtheorem{Prop}[Theorem]{PROPOSITION}
 \newtheorem{Proposition}[Theorem]{PROPOSITION}
 \newtheorem{Lemma}[Theorem]{LEMMA}

\title{Another proof of M. Kontsevich formality theorem
for ${\bbf R}^n$}
\begin{document}
\newcommand{\be}{\begin{equation}}
\newcommand{\ee}{\end{equation}}
\hsize 15cm
\vsize 20cm
\leftmargin=-4cm
\topmargin=1cm
\maketitle
\section{Introduction}
This is a draft of paper in which we announce a plan
of an alternative proof of M. Kontsevich
formality theorem \cite{K}.
The basic idea is to equip Hochschild cochains of
an associative algebra $A$ with a structure of homotopy
Gerstenhaber algebra and to prove the formality of this Gerstenhaber algebra.
See \cite{Vor}.
The author  was told about this idea by Boris Tsygan
about a year ago. 

The homological obstructions to  formality  vanish in the case
$A=SR^n; C^\infty(R^n)$. In other words, the Gerstenhaber algebras
formed by Hochschild cohomology are not deformable.  
 The operad $e_2$ governing Gerstenhaber algebras
is Koszul, therefore it has a canonical resolution
which we denote by $\he$. The Hochschild cochains
of an associative algebra have  a canonical structure
of an algebra over the operad $\bi$ (see \cite{GJ} and section \ref{descbi}), and it suffices
to construct a map $\he\to \bi$. Such a map must induce
the correct structure of Gerstenhaber algebra on the Hochschild cohomology
and  the correct Gerstenhaber bracket on Hochschild cochains. 
These conditions are formalized in Theorem \ref{mainth},
and in the sections \ref{imply}  and \ref{obstr}  
 M. Kontsevich's theorem is deduced from  Theorem \ref{mainth}.

In the section \ref{redu} we formulate Theorem \ref{third}
and show that it implies Theorem \ref{mainth}. The rest of the
paper is devoted to the proof of this theorem. 

In section \ref{bitoe} we construct a map $k:\bi\to\e$.
The existence of such a map (satisfying
the condition 1 of Theorem \ref{third})
is a direct corollary
of Etingof-Kazhdan theorem on quantization of bialgebras. 

In section \ref{stroit}
we construct the operad $\F$. If we new that the homology
operad of $\bi$ is $\e$ we could take $\F=\bi$. 
Let us outline the main steps of construction of $\F$.
Suppose we have an operad $\x$
in the category of dg-coalgebras with
counit (for example the asymmetric operad 
$\As$ from section \ref{enrich})
and a dg coalgebra with counit $W$. 
Then we can define a notion of an $\x$-algebra
on $W$. It is the same as in the case of usual operads,
but in addition we require that all structure maps be
coalgebra morphisms. Then we consider the case when $W$
is  cofree
as a graded coalgebra, that is $W=TV$, and we prove
that in this case $\x$-algebras on $W$ are determined
by a collection of operations $TV\to V$ and that
the relations between these operations are
described by a certain dg-operad, which we denote by
${\cal O}(\x)$. It is clear that $\bi$ is  constructed 
exactly in this way from $\x=\As$. 
If $\a$ is another operad of coalgebras with counit and
$f:\a\to \As$ a morphism, then we have the induced morphism
$\O(f):\O(\a)\to \O(\As)$ and our $\F$ is of the form $\O(\a)$.
Also, we prove that if $\a$ is free as a dg-operad, then
so is $\O(\a)$. 

Thus, we need an $\a$. We want $\a$ to be free as a dg operad
and we want $f$ to be a quasiisomorphism. In other words,
we are going to construct a resolution of $\As$.
It is natural to take the chain operad of Stasheff
associahedra. Since we need zero-ary operations we
have to decompose the associahedra. This is done
in section \ref{stasheff}.

In the sections \ref{cohom1}-\ref{cohom2} we prove that
the through map $\F\to\bi\to\e$ is a quasiisomorphism.
We take the filtration on $\F$ defined by
the number of internal vertices on the trees, and consider
the associated spectral sequence. This sequence 
is similar to the spectral sequence 
of manifolds with corners which are the Fulton-McPherson
compactifications of configuration spaces of $n$ points
on ${\Bbb R}^2$, see \cite{GJ}. The author believes
that $\F$ is a chain operad of the operad of configuration
spaces for a suitable decomposition. 

Finally,
in section \ref{last}
 we construct a map $\hlie\{1\}\to\F$ that insures
that this construction gives the correct
 Gerstenhaber bracket.   
   
The author would like to thank Boris Tsygan and Paul
Bressler for their help.
\section{Statement of the main theorem and why it implies the Formality
theorem}
\subsection{Generalities}

\subsubsection{Conventions} 1.We denote
by $k$ some fixed field of characteristic 0. 
 Differentials in all complexes will have degree
1 unless the opposite is stated. As usual, for any complex $V^{\bullet}$
we define the shifted complex $V[k]$ so that $V[k]^{i}=V^{k+i}$.
 The grading on the 
dual complex $V^*$ is defined
by $V^{*k}=(V^{-k})^*$.

2. For an operad $\O$ we denote by $\O(n)$
the set of $n$-ary operations. For  an operation $s\in \O(n)$, the symbol
$s(x_{i_1},\ldots,x_{i_n})$ will denote the effect of the application of the permutation $(i_1,\ldots,i_n)$ to $s$. The structure map
of insertion of an $n$-ary operation  in the $i$-th position
of an $m$-ary operation $\O(m)\otimes \O(n)\to
\O(m+n-1)$ will be denoted by $\circ_i^{m,n}$. Sometimes the upper indices
will be omitted. 

3. As usual, we denote by $\ass$, $\com$, $\lie$, $\e$ the operads governing
associative, commutative, Lie, and Gerstenhaber algebras respectively. 

\subsubsection{Shift in operads}\label{shift} For an arbitrary dg-operad $\O$, we are 
going to define an operad $\O\{m\}$, such  that a structure of
an $\O\{m\}$-algebra on a complex $V$ is equivalent to a structure of an
$\O$-algebra on $V[m]$. First, we define the operad $\com\{m\}$ as follows.
We set $\com\{m\}(n)=\Lambda_n^{\otimes m}[(n-1)m]$, where $\Lambda_n$
is the sign representation of $S_n$. Let $1_n\in \com\{m\}(n)$ be the 
generator. Set $\circ_1^{p,q}(1_p,1_q)=1_{p+q-1}$. This uniquely defines 
$\com\{m\}$. For an arbitrary dg-operad $\O$ we set $\O\{m\}=\O\otimes \com\{m\}$.
That is $\O\{m\}(n)=\O(n)\otimes \com\{m\}(n)$, and the structure maps
are defined as tensor products of the corresponding   
structure maps of  $\O$ and $\com\{m\}$.

\subsubsection{} It is well known that the operads $\ass$, $\com$, $\lie$, $\e$ are Koszul.
Therefore, they have canonical resolutions $\hass$, $\hcom$,
$\hlie$, $\he$ (see \cite{GK},
\cite{GJ}). We follow the definitions in \cite{GJ}.
We denote by $\O^\lor{}$ the Koszul dual cooperad for a Koszul operad
$\O$.
By definition from \cite{GJ},  $\O^\lor{}$ is cogenerated
by $\O^\lor{}(2)\cong \O(2)[1]$, and the corelation space
is isomorphic to 
$$
\O(2)[1]\otimes (\O(2)[1]\otimes_{S_2}k(S_3))/R[2],
$$
where $R[2]$ is the relation space of $\O$.

 A homotopy $\O$-algebra on a complex $V$ 
is the same as a differential on a cofree
$\O^\lor$-coalgebra $S_{\O^\lor}V$ cogenerated by $V$, where 
$S$ is the Schur functor. Homotopy $\O$-algebras are governed by  an operad
$\O'$ which is a semi-free operad generated by an $S$-module $\O^\lor[-1]$.
We have $\com\{n\}^\lor=\lie\{-1-n\}^*$; $\lie\{n\}^\lor=\com\{-1-n\}^*$;
$\e\{n\}^\lor=\e\{-2-n\}^*$, where $*$ denotes the linear dual cooperad.
\subsubsection{Description of $\he$}\label{hoe} This operad was described in \cite{GJ}.

A homotopy $e_2$-algebra on a graded vector
space $V$  is a differential on a cofree coalgebra $Q(V)$
over the cooperad $e_2^\lor $ cogenerated by $V$. Let $P(V)\cong S_{\e\{-1\}^*}V$ be the cofree 
$\e\{-1\}^*$-coalgebra. Then  
 $Q(V)\isom P(\V)[-1]$.
It will be more convenient for us to use $P(\V)$ rather then
$Q(V)$. 
 $P(\V)$ has two co-operations 
$P(V[1])\to P(V[1])\otimes P(V[1]) $ 
comultiplication and cobracket
with their degrees +1 and 0 respectively, and a homotopy $e_2$-
algebra on $V$ is the same as a differential on $P(\V)$.

\subsubsection{Operad $\bi$}\label{descbi}  This operad was described in \cite{GJ}.
Here we reproduce its definition in a way convenient for us. By
definition, $\bi=\b\{1\}$, where $\b$ is such an operad that a structure 
of a $\b$-algebra on a complex $V$ is equivalent to  a structure of
a dg-bialgebra on $TV$ with standard coproduct
\begin{equation}
\Delta e_1\otimes e_2\otimes\ldots\otimes e_n=\sum_{i=0}^{n}e_1\otimes\ldots
\otimes  e_i\bigotimes
e_{i+1}\otimes\ldots\otimes e_n,
\end{equation}
with  
the standard unit $1\in T^0V$, with the standard counit $\epsilon: TV\to T^0V$,
and with a differential $D:TV\to TV$ such that the restriction of it on $V\isom T^1V$ takes values in $V\isom T^1V$ and is equal to the differential
on $V$ as a complex. To specify a $\b$-algebra one has to determine the maps
$m_k:T^kV\to V,\ k\geq 2$ (corestrictions of $D$ on $V$) and
$m_{k,l}:T^kV\otimes T^lV\to V,\ k,l\geq 1 $ (corestrictions on $V$ of the product).The operad $\b$ is generated by these operations. These operations
should satisfy certain identities that provide the associativity
of the product. The condition $D^2=0$ defines the differential of $m_k$.
The condition that the product agrees with the differential
defines the differential of $m_{k,l}$. The degree of $m_{k,l}$ is 0;
the degree of $m_k$ is 1. For $\bi(n)=\b(n)\otimes\Lambda_n[n-1]$,
the degree of $\mu_{k,l}=m_{k,l}\otimes 1_{k+l}$  is $1-k-l$,
and the degree of $\mu_k=m_k\otimes 1_k$  is $2-k$.
For more details see Section \ref{functor}

\begin{Prop} $H^\bullet(\b(2))$ is generated by the class of
$m_{1,1}(x_1,x_2)-m_{1,1}(x_2,x_1)$ in $ H^0(\b(2))$ and
$1/2(m_2(x_1,x_2)-m_2(x_2,x_1))$ in $H^{1}(\b(2))$.
Respectively, $H^\bullet(\bi(2))$ is generated by the class of
 $\mu_{1,1}(x_1,x_2)+\mu_{1,1}(x_2,x_1)$ in $ H^{-1}(\bi(2))$ and
$1/2(\mu_2(x_1,x_2)+\mu_2(x_2,x_1))$ in $H^{0}(\bi(2))$. 
These generators define an isomorphism 
\begin{equation}\label{isombie}
\varkappa:H^\bullet(\bi(2))\isom \e(2)
\end{equation}
\end{Prop}\begin{pf} Direct computation
          \end{pf}
\subsection{Statement of the main theorem}
\subsubsection{} Direct check shows that the operation $m_{1,1}(x_1,x_2)-m_{1,1}(x_2,x_1) \in \b(2)$ satisfies the Jacoby identity
and its differential is zero. Therefore, we have 
maps $\lie\to \b$ and 
\begin{equation}\label{lietobi}
\lie\{1\}\to\bi.
\end{equation}
 Also we have a map $\lie\{1\}\to \e$,
corresponding to the forgetting of the commutative product. This map
defines a map of the resolutions
\begin{equation}\label{hlietohe}
\phi:\hlie\{1\}\to \he.
\end{equation} 
\subsubsection{The main theorem} \begin{theorem}\label{mainth} There exists a morphism of operads 
$q:\he\to\bi$,
such that
\begin{enumerate}
\item[1]
\begin{equation}
\begin{matrix} \he  & \stackrel{q}{\longrightarrow} & \bi\\
                \uparrow\lefteqn{\phi}&    &\uparrow\\
               \hlie\{1\}   &\longrightarrow&\lie\{1\} 
\end{matrix}
\end{equation}
\item[2]
\begin{equation}
\begin{matrix}
H^\bullet(\he) & \stackrel{q^*}{\longrightarrow} & H^{\bullet}(\bi(2))\\
\sim\downarrow   &       \swarrow   \varkappa       &    \\
\e(2)              &                         &    \\ 
\end{matrix}
\end{equation}
\end{enumerate}
\end{theorem}
\subsubsection{Why does this theorem imply formality?}
\label{imply} Let $A=SV$
be the symmetric algebra for a graded vector space $V$
(the cases $A=C^{\infty}V$ or $A=k[[V]]$ are treated in the same way). Then $C^{\bullet}(A,A)$
is naturally a $\bi$-algebra (see \cite{GJ}). Hence, via $q$, it is also an $\he$-algebra.
The cohomology $HH^\bullet(A,A)\cong SV\otimes \wedge V$ is an $\e$-algebra.
Condition 2 in theorem \ref{mainth} assures that this is the Schouten algebra.
In section \ref{obstr} we  prove that there are no obstructions to the formality
of $C^\bullet(A,A)$
as an $\he$-algebra. Therefore $C^\bullet(A,A)$ is formal as an $\he$-algebra.

 Denote $W=C^\bullet(A,A)$, $S=HH^\bullet(A,A)$. Let $P(\bullet)$ denote
the same coalgebra as in \ref{hoe}. Then the structures of
$\he$-algebras on $W$ and $S$ define differentials on $P(W[1])$ and $P(S[1])$.
The formality of $W$ implies that there exists a map of $e_2^\lor$- coalgebras
$f:P(S[1])\to P(W[1])$, such that its restriction on   $S[1]$ gives 
a quasiisomorphism $S[1]\to W[1]$. 

Let $X,Y$ be dg-spaces. The cocommutative coproduct on 
$P(X)$ defines a structure of $\com\{-1\}$-coalgebra on $P(X)$.
We have a canonical inclusion
of $\com\{-1\}$-coalgebras $i:S^{\geq 1}(X[1])[-1]\to  P(X)$
corresponding to the 'co-forgetting' of the cocommutator on $P(X)$.
One can easily prove that under any coalgebra morphism $P(X)\to P(Y)$
$i(S^{\geq 1}(X[1])[-1])$ goes to $i(S^{\geq 1}(Y[1])[-1]$ and that 
$i(S^{\geq 1}(X[1])[-1]$
is preserved by any coderivation of $P(X)$ (as a subspace, not pointwise).
In particular, this implies that any $\he\{-1\}$-algebra on $X$ 
(which is the same as a differential on $P(X)$) defines
a $\hlie$-algebra on $X$ (which is a differential 
on $S^{\geq 1}(X[1])[-1]$). The corresponding map $ \hlie\to\he\{-1\}$
is induced by the map \eqnref{hlietohe}.

 We have the following diagram of maps of $\com\{-1\}$-coalgebras.
\begin{equation}\label{diaglie}
\begin{matrix}
P(S[1]) & \longrightarrow & P(W[1])\\
\lefteqn{i}\uparrow &    &\lefteqn{ i}\uparrow\\
S^{\geq 1}(S[2])[-1]&\longrightarrow& S^{\geq 1}(W[2])[-1]
\end{matrix}
\end{equation}
Differentials on $S^{\geq 1}(S[2])[-1]$ and $S^{\geq 1}(W[2])[-1]$
are induced by the ones on $P(S[1])$  and  $P(W[1])$. 
The bottom map is induced by the top one. Note that
$S^{\geq 1}(S[2])[-1]$ is nothing else but the chain complex of the Schouten
algebra as a Lie algebra. Condition 1 in theorem \ref{mainth} implies that
$S^{\geq 1}(W[2])[-1]$ is the chain complex of the Gerstenhaber algebra on $W$.
The bottom map of the diagram \eqnref{diaglie} gives the quasiisomorphism
of Schouten algebra and Gerstenhaber algebra since its restriction
onto $S[1]$ coincides with the restriction on $S[1]$ of the top map.

\section{Obstructions to the formality}\label{obstr}
\subsection{Additional gradings on $P(V[1])$} 

According to the description in the section \ref{hoe}, an $\he$-algebra on $V$
is the same as a differential of
$P(\V)$. The definition of $P(\V)$ implies that
$P(\V)\cong S((T(V[1])/shuffles)[1])[-1]$ as a vector space.
Introduce  two additional gradings $\gr_2$ and $\gr_3$
on $P(\V)$: by the number
of cobrackets  and by the number of comultiplications.
 The grading
$\gr_2$ is first defined to be $k-1$ on $(T^k(V[1])/shuffles)[1]$
and is extended additively on $P(\V)$. 
We define $\gr_3|_{S^k((T(V[1])/shuffles)[1])[-1]}=k-1$.
Thus, $P(V[1])$ is  a 3-graded vector space.
the first grading $\gr_1$ is the grading of $P(\V)$ as
an $(e_2\{-1\}^\lor)$-coalgebra, the second grading
is the number of cocommutators, and the third
grading is the number of comultiplications.
Comultiplication has degree $(1, 0, -1)$; cobracket has degree $(0,-1,0)$.

\subsection{Complex of coderivations}

  Coderivations of $P(\V)$ form
a Lie super-algebra. Any coderivation of $P(\V)$  is uniquely defined by its corestriction onto $\V$. In other words
$$
{\rm Coderivations\ of}\ P(\V)\isom {\rm Hom}_k(P(\V),\V).
$$ Set  $\gr_2(\V)=\gr_3(\V)=0$; $\gr_1(\V)=\gr(V)-1$, where $\gr$ is the original grading on $V$. Then ${\rm Hom}_k(P(\V),\V)$ becomes a trigraded space and this grading of any  derivation coincides with grading of this derivation as
an element of  ${\rm Hom}_k(P(\V),P(\V))$. An $e_2-$algebra on $V$ is the same as a differential of $P(\V)$ centered at the gradings $(1, -1, 0)$ (multiplication) and $(1, 0, -1)$ 
(bracket). Denote the 
$(1, -1, 0)$-part of this differential by $d_m$, and the $(1, 0,-1)$-part by $d_{br}$.
Both of these parts are also differentials (since they also determine $e_2$-algebras).
 Therefore, $(({\rm Hom}_k(P(\V),\V)^{\gr_2, \gr_3}, [d_m,\cdot],[d_{br},\cdot])$ is a bicomplex. Here the brackets denote
the commutator in the Lie algebra of coderivations.
Cohomology groups of this bicomplex  in degrees less than 0
are obstructions to formality
of a homotopy  $e_2$-algebra 
whose cohomology algebra is $V$.
 We will prove that there is no such cohomology
in the case when $V$ is the Schouten algebra.
Since this bicomplex is concentrated in  negative degrees, the spectral sequence of it
converges. 
\subsection {$E_1^{\bullet,\bullet}$}

We have to compute the cohomology with respect to $d_m$. Let us describe the action of this differential restricted to the part of ${\rm Hom}_k(P(\V),\V)$ with grading $\gr_3=1-k$ which is isomorphic
to
$$
{\rm Hom}_k(S^k(T(\V)/shuffles[1]))[-1],V[1]),
$$
where all the shifts here and below are made with respect to $\gr_1$.
Let $({\rm Harr}(V[1]),b)$ be the standard Harrison complex  of $V[1]$ as a $V$-module, and $V$ as a cocommutative coalgebra:
$$
{\rm Harr}(V[1])=((T(V[1])/shuffles)\otimes V[1], b),
$$
 where $b$ is induced by 
Hochschild differential.
We define $\gr_1$ and  $\gr_2$ on
 ${\rm Harr}(V[1])$
by setting $\gr_1|_{V[1]}$ to be the original grading on $V[1]$, $\gr_2|_V=0$.
We define $\gr_1$  on $T(V[1])/shuffles$ as the induced grading from $V$
and we set $\gr_2|_{T^k(V[1])/shuffles}=k-1$.
 Then ${\rm Harr}(V[1])$  is a complex with respect to each of
these gradings. We have
$$
{\rm Hom}_k(S^k(T(\V)/shuffles[1])[-1],V[1])\cong
{\rm Hom}_{V^{\otimes k}}(S^k({\rm Harr}(V[1]))[-1],V[1]).
$$
\subsection{ Schouten algebra.}
Let $V$ be the  Schouten algebra. Then as a commutative algebra $V=S(W)$ for a certain finite dimensional
graded space $W$. In this case, it is well known that
${\rm Harr} (V[1])$ is quasiisomorphic to 
$W[1]\otimes V[1]$ with $\gr_1$ induced by the gradings
on $V$ and $W$, and $\gr_2=0$. Therefore,
our complex is quasiisomorphic to 
$$
{\rm Hom}_{V^{\otimes k}}(S^k(W[1]\otimes V[1])[-1],V[1])= V[2]\otimes 
S^k((W[2])^*),
$$
This cohomology has grading $(\gr_2=0,\gr_3=1-k)$ and an element $\phi=f\otimes w_1\otimes\ldots\otimes w_k$,\quad
$f\in V; w_i\in W^*$, is presented by the following map $m_\phi:P(\V)\to \V$:
$m_\phi$ is non zero only on  $P(\V)^{(\gr_2=0,\gr_3=1-k)}\isom S^k(V[2])[-1]$, 
\begin{equation}\label{represent}
m_\phi(f_1,f_2,\ldots  ,f_k)=fi_{w_1}f_1\cdot\ldots\cdot i_{w_k}f_k(-1)^{\sum \gr w_i(\gr f_1+\ldots \gr f_{i-1})}
\end{equation}

Thus, $E_1^{0,k}=V[2]\otimes S^{1-k}(W[2]^*)$.
\subsection {$E_2^{\bullet,\bullet}$}

Differential $d_{br}$ restricted to the zeroth line $\gr_1=0$ coincides with the Lie differential with respect to the bracket.  The space $W$ is equal to $U\oplus U^*[-1]$, where $U=R^n$.
Therefore $W^*=W[1]$.
Hence, our $E_1^{\bullet,\bullet}$ is $ S W[2]\otimes S^{\geq 1}(W [-1]) $
 and it is easy to see that
the differential $d_2$ induced by $d_{br}$ is just the de Rham differential sending $W$ identically to $W[-1]$. Since
 the zeroth power in the right multiple is truncated, 
this complex  has cohomology
of grading (0,0) and no other cohomology. This cohomology means that we can deform 
the differential in the Schouten algebra to be the bracket with some element of it.
This elements do not create an obstruction to formality.

\section{First Reductions}\label{redu} Let us state a theorem which implies theorem
\ref{mainth}
\begin{theorem}\label{second} There exists an operad $\F$, a quasiisomorphism
$qis:\F\to\e$, a map $r:\F\to\bi$, and a map $s:\hlie\to \F$, such that
the following diagrams are commutative.
\begin{enumerate}
\item[1]
\begin{equation}\label{pervoe}
\begin{matrix}
\bi & \stackrel{r}{\longleftarrow} & \F     & \stackrel{qis}{\to}  & \e\\
\uparrow&                      &\uparrow&                     & \uparrow\\
 \lie\{1\}&\longleftarrow             &\hlie\{1\}&     \stackrel{\eqnref{hlietohe}}{\to}                                                              &\he
\end{matrix}
\end{equation}
\item[2]

\begin{equation}
\begin{matrix}
H^\bullet(\F(2)) & \stackrel{qis_*}{\longrightarrow} & \e(2)\\
\downarrow \lefteqn{r_*}   &\nearrow\lefteqn{\sim}           & \\
H^\bullet(\bi(2))&                       &
\end{matrix}
\end{equation}
\end{enumerate}
\end{theorem}
\subsection {Theorem \ref{second} implies theorem \ref{mainth}}
We are going to use the structure of closed model category
on the category of dg k-operads, where $k$ is a field of characteristic 0. 
\begin{Lemma} The map \eqnref{hlietohe} is a cofibration.
\end{Lemma}
\begin{pf} It follows from the fact that $\he$ can be
obtained from $\hlie\{1\}$ by adding free variables and by killing
the cycles.
\end{pf}

Now, consider the diagram (\ref{pervoe}).
Using RLP, we have a map $\phi:\he\to \F$, such that the composition
$q=r\circ\phi$ satisfies the conditions of theorem \ref{mainth}.
\subsection{One more reduction}
\begin{theorem}\label{third}
\begin{enumerate}
\item[1] There exists a morphism $k:\bi\to e_2$, such that $k_*:H^\bullet(\bi(2))\to H^\bullet(\e(2))$ coincides with the map \eqnref{isombie}.
\item[2] There exists an operad $\F$ and morphisms $s:\hlie\{1\}\to\F$
and $t:\F\to \bi$, such that the diagram
\begin{equation}
\begin{matrix}
       &                   &\e\\
       &\nearrow \lefteqn{r}         &\uparrow\lefteqn{k}\\
\F     &\stackrel{t}{\longrightarrow}  & \bi\\
s\uparrow&                 &\uparrow\lefteqn{\eqnref{lietobi}}\\
\hlie\{1\}&\longrightarrow             &\lie\{1\}
\end{matrix}
\end{equation}
is commutative
and 
\item[3] the map $r=k\circ t$ is a quasiisomorphism.

\end{enumerate}
\end{theorem}  
It is clear that this theorem implies theorem \ref{second}. 

\section{Proof of theorem \ref{third} 1} \label{bitoe}
Any map $\bi\to \e$ defines  a functor from
the category of $\e$-coalgebras to the category of $\b$-coalgebras.
We will construct such a functor, and the map $\bi\to \e$ will be defined
as a unique map producing this functor.
\subsection{Category of $\bi$-coalgebras}. Let $V$ be a graded  vector space.
Set $\pt V[1]=\prod\limits_{k=0}^{\infty}T^kV[1]$. Endow this space
with the p-adic topology in which the base of open neighborhoods
of 0 is formed by the subspaces $\prod\limits_{k=l}^{\infty}T^kV[1]$.
The space $\pt V[1]$ has a natural structure of a topological algebra with
 unit. Indeed, the multiplication is extended by continuity from
the one on $TV[1]$, and the unit is $1\in T^0V[1]$. Define $\epsilon:\pt V[1]\to k$ as the projection on $T^0V[1]$.
\begin{Proposition} The category of $\bi$-coalgebras is equivalent to the 
category whose objects are topological dg-bialgebras isomorphic to $\pt V[1]$
as algebras with unit and having $\epsilon$ as a counit.  Morphisms
between $\pt V[1]$ and $\pt W[1]$ are maps $V\to W$ which induce
a morphism of topological dg bialgebras with unit and counit.
\end{Proposition} 
\subsection{Construction of the Functor}
\subsubsection{Etingof-Kazhdan theorem} We are going to use a result from
\cite{KE}. Let us reformulate it. Let $\frak a$ be a Lie bialgebra with commutator $[,]$ and cocommutator $\delta$. Let $U({\frak a})$ be the enveloping algebra
of $\frak a$ as a Lie algebra.
\begin{theorem}(Etingof-Kazhdan)\label{EtKa} There exist
k[[h]]-linear  maps $m:U({\frak a})\otimes U({\frak a})[[h]]\to U({\frak a})[[h]]$, 
$\Delta:U({\frak a})[[h]]\to U({\frak a})\otimes U({\frak a})[[h]]$ such that 
\begin{enumerate}
\item[1] $(U({\frak a})[[h]], m, \Delta)$ with the standard unit $1\in U({\frak a})[[h]]$and counit $\epsilon:U({\frak a})[[h]]\to k[[h]]$ is a bialgebra 
with unit and counit over $k[[h]]$.
\item[2] Via  the PBW identification $U({\frak a})\sim S({\frak a})$, the components of $m$ $m^r_{p,q,l}:S^p{\frak a}\otimes S^q{\frak a}\to S^r({\frak a})h^l$
and of $\Delta$ $\Delta^{p,q}_{r,l}:S^r({\frak a})\to S^p({\frak a})\otimes S^q({\frak b})h^l$
are obtained from the operations $[,]$ and $h\delta$
via the acyclic calculus over $\bbf Q$ (without using $h$).
\item[3] The multiplication coincides with the usual multiplication
on $U({\frak a})[[h]]$ up to $O(h)$.
\item[4] $\Delta$ is the standard coproduct on $U({\frak a})[[h]]$
up to $O(h)$
\item[5] On ${\frak a}\subset U({\frak a})$
we have $\Delta{\frak a}-\Delta^{op}{\frak a}=h\delta a+O(h^2)$
\end{enumerate}
\end{theorem} 
 One sees that this theorem is also applicable to Lie dg bialgebras.
In this case we obtain a structure of a dg-bialgebra with unit and counit
on $U({\frak a})[[h]]$, and the differential on $U({\frak a})[[h]]$ is induced
from the one on $\frak a$.
\subsubsection{A dg Lie bialgebra $\f V[1]$} 
We  construct a dg Lie bialgebra to which we will apply the Etingof-Kazhdan
theorem.

Let $V$ be a dg $\e$-coalgebra. This means that we have a cocommutative coproduct $\Delta: V\to V\otimes V$, a cocommutator $\delta: V[1]\to V[1]\otimes V[1]$
and a differential $d:V\to V$. In particular, $V[1]$ is a dg Lie coalgebra.
Let $\f V[1]$ be the free Lie algebra generated by $V[1]$. Define a cocommutator
$\bar\delta$ on it as a unique map $\bar\delta:\f V[1]\to \f V[1]\otimes \f V[1]$ such that 
\begin{enumerate}
\item[1] $\bar\delta([x,y])=[\bar\delta x,y]+[x,\bar\delta y]$;
\item[2] $\bar\delta|_{V[1]}=\delta.$
\end{enumerate}
One can check
 that $\bar\delta$ turns $\f V[1]$ into a Lie bialgebra.
This bialgebra was introduced in \cite{Dr}. The coproduct $\Delta$ and
the differential $d$ define a bar-differential $b$ on the Harrison complex
of a dg cocommutative coalgebra $V$. This complex, as a vector space, is
isomorphic to $\f V[1]$. One can check that $b$ is compatible with 
the structure of a Lie bialgebra on $\f V[1]$. Thus, $\f V[1]$ has a natural
structure of a dg Lie bialgebra. 

\subsubsection{Application of the Etingof-Kazhdan Theorem}
The Etingof-Kazhdan
theorem gives a structure of a dg bialgebra on $U(\f V[1])[[h]]\cong TV[1][[h]]$ with respect to
the grading $\rm gr$ induced form the one on $V$ and such that
${\rm gr}\;h=0$. The differential $\bar b$  is induced from the bar-differential 
$b$ on $\f V[1]$ and is nothing else but the bar differential of $V$ as 
an associative dg coalgebra.      

Note that $\f V[1]$ has another grading $| \ |$ such that
$|[e_1,[e_2,\ldots, e_k]\ldots]|=k$. In this grading $|[,]|=0$
and $|\delta|=1$. The differential can be split into two parts $b=b_0+b_1$
with gradings 0 and 1 respectively ($b_0$ corresponds to the differential
on $V$ and $b_1$ corresponds to the commutative coproduct $\Delta$).
This grading naturally extends to $U(\f V[1])\cong TV[1]$ in such a way that
$|e_1\otimes\ldots\otimes e_k|=k$. Set $|h|=-1$. Then the gradings of
the product and the coproduct on $TV[1][[h]]$ are equal to 0. This allows
one to set $h=1$ and to extend the product and the coproduct on
$\pt V[1]$ by continuity. The bar-differential $\bar  b$ can be
also extended to $\pt V[1]$ by continuity. Therefore, we have obtained
a structure of a topological dg-bialgebra on $\pt V[1]$. The only
thing that prevents it from being a $\bi$-algebra is the fact that 
the product $m$ on $\pt V[1]$ is non-standard one. The Etingof-Kazhdan
theorem implies the following:
\begin{enumerate} 
\item[1] $m$ has a unit $1\in \pt V[1]$;
\item[2] the counit map $\epsilon:\pt V[1]\to k$ is a morphism of algebras;
\item[3] for $x\in T^kV[1]$, $y\in T^lV[1]$, we have $m(x,y)=x\otimes y+z$,
where $z\in T^{>k+l}V[1]\subset \pt V[1]$.
\end{enumerate}
The last statement follows from theorem \ref{EtKa}, 3.

Any product $m$ satisfying these conditions is isomorphic to the standard
one. Denote $x*y=m(x,y)$.  The isomorphism $\phi$ is given by the formula 
$\phi(e_1\otimes\ldots\otimes e_k)=e_1*\ldots*e_k,\ k\geq 1 $, $e_i\in V[1]$;
$\phi(1)=1$. Since $\phi(e_1\otimes\ldots\otimes e_k)\in \pt^{\geq k}V[1]$,
 $\phi$ is well defined on $\pt V[1]$. Also, for  an $x$ such that $|x|=k$, 
we have $\phi(x)=x+y$, $y\in \pt^{>k}V[1]$.
Therefore, $\phi^{-1}$ is a continuous map $\pt V[1]\to \pt V[1]$.
We can conjugate all structure maps by $\phi$ and we obtain a dg-bialgebra
$(\pt V[1],\phi_*m,\phi_*\Delta, \phi_*d)$ which gives a structure
of a $\bi$-algebra on $V$. Thus, we have constructed a functor
providing a map $\bi\to \e$. The quasi-classical limits from theorem \ref{EtKa}
4,5 allow one to easily  check  the last condition of theorem \ref{third} 1.  

\section{Proof of Theorem \ref{third},2: construction of $\F$}
\label{stroit}  
We will start with a description of a functor from 
a certain sub-category of the category of operads
of dg  coalgebras with counit to the category of dg operads such that both
$\bi$ and $\F$  are in the image of this functor. 
\subsection{Description of the functor}\label{functor}
\subsubsection{Operads of coalgebras and algebras over them}\label{cat} Let $\x$ be an asymmetric operad in the symmetric monoidal category
of associative coalgebras with counit such that $\x(0)=\x(1)=k$. Let 
$\cx$ be the category of such operads  whose morphisms act identically on
the spaces of unary and 0-ary operations. Let $W$ be a dg-coalgebra with
counit. Let $e\in W$ be a marked element.
\begin{Definition}\label{def} An $\x$-algebra with unit $e$ on $W$ is a collection of maps
$\phi_n:\x(n)\otimes W^{\otimes n}\to W$ such that
\begin{enumerate}
\item[1] Each $\phi_n$ is a morphism of coalgebras with counit.
\item[2] Collection of $\phi_n$ is a representation
of $\x$ as a dg operad on a dg space $W$.
\item[3] $\phi_0(1,1)=e$; $\phi_1(1,x)=x$
\end{enumerate}
\end{Definition}
Let $V$ be a complex.
 Let $TV$ be a coalgebra with  coproduct, counit, unit, and a differential 
described in section \ref{descbi}. Note that the differential is not defined
uniquely.
\begin{Definition}
An $\x\b$-algebra on $V$ is an $\x$-algebra with unit on $TV$.
\end{Definition}   
\subsubsection{Description of $\x\b$-algebras} To define an $\x\b$-algebra
on $V$ one has to specify
a differential $D:TV\to TV$
satisfying the conditions from section \ref{descbi} and   maps $\phi_n:\x(n)\otimes TV^{\otimes n}\to TV$
satisfying \ref{def}. Let us decompose $D$ and $\phi_n$ into their components
$D_k^l:T^kV\to T^lV$ and
$\phi^r_{k_1,\ldots,k_n}:T^{k_1}V\otimes\ldots\otimes T^{k_n}V\to T^rV$.
\begin{Lemma} We have
\begin{enumerate}
\item[1] $\phi^r_{k_1\ldots k_n}=0$ for $r> k_1+\ldots+k_n$;
\item[2] $D_k^l=0$ for $l>k$; $D_k^k$ is the differential induced
by the differential $d$ on $V$.
\end{enumerate}
\end{Lemma}
\begin{pf}
1. We will use double induction with respect to $k_1+\ldots+k_n$
and $n$.

a) $n=0,1$. In this cases the statement follows from the Definition \ref{def},
3;

b) Suppose that the statement has been proven for all $n\leq N$, $N\geq 2$, and
any sets of $k_1\ldots k_n$ and for $n=N$ under condition that
$k_1+\ldots+k_n<M$. Let us prove the statement for $n=N$, $k_1+\ldots+k_n=M$.

\noindent i) $M=0,1$. In this case one of the numbers $k_1\ldots k_n$ must be equal to 0.
Let $k_i=0$. Let $x_i\in T^{k_i}V$ and $a\in \x(n)$. Then 
$$
\phi_{k_1\ldots k_n}^r(a,
x_1,\ldots, x_n)=x_i\phi_{k_1\ldots \hat{k_i}\ldots k_n}^r
(\circ^{n,0}_i(a,1), x_1,\ldots \hat{x_i}\ldots, x_n)
=0
$$
 if $r>k_1+\ldots+k_n$ by the induction assumption.

\noindent ii) $M>1$.   Let $x_i\in T^{k_i}V$ and $a\in \x(n)$.
Let $\Delta x_i=\sum \Delta^{p,k_i-p}x_i$, where $\Delta^{p,k_i-p}x_i\in
T^pV\otimes T^{k_i-p}V$. Then  
\begin{multline}
\Delta \phi_n(a, x_1,\ldots, x_n)-
1\otimes\phi_n(a, x_1,\ldots, x_n)-
\phi_n(a, x_1,\ldots, x_n)\otimes 1\\
=\sum \phi_{p_1\ldots p_n}^i\otimes\phi_{k_1-p_1\ldots k_n-p_n}^j
(\Delta a, \Delta^{p_1,k_1-p_1}x_1,\ldots,
\Delta^{p_n,k_n-p_n}x_n),
\end{multline}
where the sum is taken over all $i,j,p_1,\ldots,p_n$ such that
$0\leq p_i\leq k_i$, not all $p_i$ are equal to 0, and
 not all $p_i$ are equal to $k_i$.
By the induction assumption the right hand side belongs to
$\bigoplus\limits_{r+s\leq k_1+\ldots+k_n} T^rV\otimes T^sV$. Therefore,
    $\phi_n(a, x_1,\ldots, x_n)\in 
T^{\leq  k_1+\ldots+k_n}V$. This proves the first statement of Lemma.
The second statement can be proved similarly.
\end{pf}
\begin{Proposition}\label{struct}
A structure of an $\x\b$-algebra on a complex $(V,d)$ is the same as 
a collection of maps $\phi^r_{k_1\ldots k_n}:\x(n)\otimes T^{k_1}V\otimes\ldots\otimes
T^{k_n}V\to T^rV$, $r\leq k_1+\ldots+k_n$, $k_i\geq 0$,
$n\geq 0$ of degree 0 and
$D_l^m: T^lV\to T^mV$, $l>m>0$ of degree 1 satisfying the following identities.
Let $x_i\in T^{k_i}V$ and $a\in \x(n)$.
\begin{enumerate}
\item[1] 
For any $i$ and $r$,
\begin{multline}\label{cond1}
\phi^r_{k_1\ldots k_n}(a, x_1,\ldots, x_n)\\
= \sum\limits_{p_1\ldots p_n }\phi^i_{p_1\ldots p_n}\otimes\phi^{r-i}_{k_1-p_1\ldots k_n-p_n}(\Delta a, \Delta^{p_1,k_1-p_1} x_1,\ldots,
\Delta^{k_n,p_n-k_n}x_n).
\end{multline}
\begin{multline}
D_l^{p+l-r}(x_l)=(D^p_r\otimes {\rm Id})(\Delta^{r,l-r} x)\pm 
({\rm Id}\otimes D^{l-r}_{l-p})(\Delta^{p,l-p}x).
\end{multline}
 These identities express the fact that $\phi_n$ are coalgebra morphisms
and that $D$ is a derivation of the coalgebra $TV$. 
\item[2] Let $b\in \x(i)$
and $c\in \x(n-i)$. Then

\begin{multline}
 \phi^r_{k_1\ldots k_n}(\circ^{n-i,i}_j(c, b), x_1,\ldots,x_n)\\
=
\sum\limits_{s\leq k_{j}+\ldots+ k_{j+i-1}}\phi^{n-i}_{k_1,\ldots, k_{j-1},s,k_{j+i},\ldots, k_n}(c, x_1,\ldots, x_{j-1},\phi^s_{k_{j}\ldots k_{j+i-1}}(b, 
x_j,\ldots, x_{j+i-1}), x_{j+i},\ldots, x_n)
\end{multline}
 
\item[3]
The map $\phi_0$ maps $\x(0)\otimes T^0V$ identically to $T^0V$.
The maps $\phi^i_j$ are identical  if $i=j$ and equal to zero otherwise.
We have $\phi^0_{k_1\ldots k_n}(a, x_1,\ldots,x_n)=0$
if not all of $k_i$ are zeros. Otherwise,  $\phi^0_{0\ldots 0}(a, x_1,\ldots,
  x_n)=x_1\cdot\ldots\cdot x_n\varepsilon(a)$.
\item[4] $dD^k_lx+D^k_ldx+\sum\limits_{l<r<k}D^k_rD^r_lx=0$
\item[5]
\begin{multline}\label{differen}
 d\phi^r_{k_1\ldots k_n}(a, x_1,\ldots, x_n)-  
\phi^r_{k_1\ldots k_n}(da, x_1,\ldots, x_n)\\
-\sum(-1)^{|a|+|x_1|+\ldots+|x_{l-1}|}\phi^r_{k_1\ldots k_n}(a, x_1,\ldots,
 dx_l,\ldots, x_n)\\
=-\sum\limits_s D^r_s\phi^s_{k_1\ldots k_n}(a, x_1,\ldots,x_n)
 +\sum(-1)^{|a|+|x_1|+\ldots+|x_{l-1}|}\phi^r_{k_1\ldots s\ldots k_n}(a, x_1,\ldots, D_l^s x_l,\ldots, x_n)\\
\end{multline}
\end{enumerate}
\end{Proposition}
\begin{pf} Clear \end{pf}
\begin{Proposition}.\label{PROP} $\x\b$-algebras are governed by a certain dg  PROP
$P(\x)$ generated
by
\begin{enumerate}
\item[1] the operations $\phi(x)^r_{k_1\ldots k_n}$, $r\leq k_1+\ldots+k_n$, $k_i\geq 0$,
$n\geq 0$, where $x\in \x(n)$ is a homogeneous element.
 Degree of such an operation is equal to $|x|$. 
\item[2]
$D_l^m$, $l>m>0$ of degree 1.
\end{enumerate}
The relations are expressed in Proposition \ref{struct} 1-3. 
The differential is defined in Proposition \ref{struct} 4,5.
\end{Proposition}
\begin{pf} The only thing that has to be checked here is that 
 the differential  respects the relations in $P(\x)$
 and that its square is equal to zero 
modulo the relations in $P(\x)$ which is obvious.
\end{pf}
\subsubsection{$P(\x)$ is generated by an operad}
\begin{Lemma}\label{odin} In Proposition \ref{struct} 2 it suffices to set $r=1$.
\end{Lemma}
\begin{pf} Conditions 1,3 in Proposition \ref{struct} imply that 
after  summation of \eqnref{cond1} over $r$, we will get on both sides
morphisms of coalgebras with counit $\x(n)\otimes TV^{\otimes n}\to TV$. Any
such a morphism is uniquely defined by its corestriction on
the cogenerators $V\in TV$.
\end{pf}

 \begin{Proposition} The PROP $P(\x)$ is freely generated by an operad
$\O(\x)$ formed
by its $(n,1)$-ary operations. (The operad $\O(\x)$ does not have to be free).
\end{Proposition}
\begin{pf}
The operad $\O(\x)$ is generated by $\phi(x)^1_{k_1\ldots k_n}$ with
$k_1+\ldots+k_n\geq 1$ and $D^1_l$, $l>1$. Conditions 1,3 in Proposition
\ref{struct} allow one to express the other operations in the PROP
in terms of these ones. They impose no restrictions on the operad $\O(\x)$.
Due to the Lemma \ref{odin}, condition 2 in Proposition \ref{struct}
can be expressed in terms of $\O(\x)$. 
\end{pf}     
\subsubsection{The case when $\x$ is semi-free as a dg operad }
\begin{Proposition}\label{liberte} Assume $\x(n)$, $n\geq 1$ is semi-free as a dg operad and is freely
generated by subspaces $\v(n)\subset \x(n)$. Then the PROP $P(\x)$
and the operad  $\O(\x)$
are also semi-free and are  freely generated by $\phi(v)^1_{k_1\ldots k_n}$,
$v\in\v(n)$, $k_1,\ldots, k_n\geq 1$, $n\geq 2$ and  $D^1_k$, $k>1$.
\end{Proposition}
\begin{pf} All relations in $\O(\x)$ are expressed in Proposition 
\ref{struct}. Now we can uniquely restore all
operations in $\o(\x)$ by induction on the number of arguments
 of an operation in $\o(\x)$.
\end{pf}
\subsubsection{Functor $\O()$.} Obviously, $\O()$ is a functor  form the category
$\cx$ described in  \ref{cat} to the category of dg operads.   \subsubsection{An asymmetric operad $\as$ and the operad $\b$.}
\label{enrich} An asymmetric operad $\as$ governing associative algebras can be enriched
in the following way. First we can turn it into an operad with 0-ary 
operations by setting $\as(0)=k$ and by setting the structure maps
$\circ^{p,q}_r:\as(p)\otimes\as(q)\to \as(p+q-1)$, $p,q,p+q-1\geq 0$,
$1\leq r\leq p$ to be  $\circ^{p,q}_r(1\otimes 1)=1$. Define coalgebra
structure on each of $\as(k)$ by $\Delta 1=1\otimes 1$. This turns $\as$
into an operad with 0-ary operations of coalgebras with counit.   It is clear that $\b=\O(\as)$.

\subsection{Operad $\a$}
 We will construct an operad of dg coalgebras with counit $\a$ with the 
following properties:
\begin{enumerate}
\item[1] $\a$ is free as an operad of graded vector spaces;
\item[2] the counit map $\varepsilon:\a\to \as$ is a quasiisomorphism 
of operads of dg coalgebras with counit.
\end{enumerate}
Then we will define $\F$ as $\O(\a)\{1\}$. The map $\O(\varepsilon)\{1\}$
will give us a map $\F\to \bi$.
 We will construct $\a$ as an operad
of chain complexes of Stasheff associahedra with respect to a suitable
polyhedral decomposition. 
 
\subsection{Stasheff topological operad as an operad with 0-ary operations}\label{stasheff}

Let $\k$  be the Stasheff topological asymmetric operad 
of associahedra. Set $\k(1)=\k(0)={\rm pt}$ and define on $\k$ a structure of an 
asymmetric operad with 0-ary operations.
For this we need a polyhedral decomposition of $\k$. The structure maps will
be defined to be piecewise linear with respect to this decomposition. Here and below 'piecewise linear map
with respect to decomposition' means 'piecewise-linear map,
mapping an element of the decomposition {\em onto} an element of a
decomposition'.

\subsubsection{Polyhedral decomposition of $\k$} Let us construct it by
induction.

\noindent1) $n\leq 2$. $\k(n)={\rm pt}$ and is decomposed trivially.

\noindent2) Suppose that all $\k(i)$, $i\leq n$ have been decomposed 
so that all operadic maps $\circ^{p,q}_l:\k(p)\times \k(q)\to \k(p+q-1)$,
$1\leq p,q, p+q-1\leq n$ are piecewise linear with respect to this decomposition.
 Then

\noindent i) decompose the boundary $\partial\k(n+1)$ which is 
  the union of $n-2$-dimensional faces of the form $\k(i)\times\k(j)$,
$i+j=n+2$, $i,j\geq 2$. Decompose each such a face  
as a product of already decomposed spaces $\k(i)$ and $\k(j)$. Let us check
that these decompositions agree on  intersections of the faces. 
Any such an intersection is a triple product $\k(p)\times\k(q)\times\k(r)$,
$p+q+r=n+3$, $p,q,r\geq 2$. The inclusions of these intersections into 
the faces look like
\begin{equation}
\begin{matrix}
(\k(p)\times\k(q))\times \k(r)&
            \stackrel{\circ^{p,q}\times {\rm Id}}{\hookrightarrow}
                                                      & \k(p+q-1)\times \k(r)\\
\downarrow \lefteqn{\sim}                &                      &\\  
  \k(p)\times(\k(q)\times\k(r)) &
                 \stackrel{{\rm Id}\times \circ^{q,r}}{\hookrightarrow}
                                                      & \k(p)\times\k(q+r-1)
\end{matrix}
\end{equation}
By the induction assumption, each of these inclusions is a piecewise 
linear map, therefore, the decompositions do agree.

ii) Let $O_{n+1}$ be the center of $\k(n+1)$. Decompose $\k(n+1)$ as the union
of cones with vertex $O_{n+1}$ and bases  the polyhedra of the decomposition of $\partial \k(n+1)$. 
It is clear that the operadic maps $\circ^{p,q}_l:\k(p)\times \k(q)\to \k(p+q-1)$,
$1\leq p,q, p+q-1\leq n+1$ are piecewise linear with respect to such a decomposition. Therefore the induction assumption is true.
\subsubsection{0-ary operations} We need to define maps $\circ^{i,0}_j:\k(i)\times \k(0)\to \k(i-1)$, $1\leq j\leq i$. Again, we will do it by induction.

1) Define $\circ^{\bullet,0}_\bullet$ on $\k(2)$ and $\k(1)$ as identical maps
${\rm pt}\to {\rm pt}$.

2) Suppose that $\circ^{i,0}_j$ have been defined for $i<n$ so that all
$\circ^{p,q}_j:\k(p)\times \k(q)\to \k(p+q-1) $, $0\leq p,q,p+q-1<n$,
satisfy the operadic identities whenever all terms in them are $\k(i)$,
$0\leq i<n$. Also suppose that each $\circ^{i,0}_j:\k(i)\times\k(0)\to \k(i-1) $ 
maps the center $O_i\times{\rm pt} \in \k(i)\times \k(0)$ to the center
$O_{i-1}$ of $\k(i-1)$.

i) Construct $\circ^{N,0}_p$ on the boundary of $\k(n)$ 
as follows. Each face on the boundary is $F=\circ^{l,m}_j(\k(l)\times \k(m))$. Define
$\circ^{n,0}_p$ on F to be
$\k(l)\times(\k(m)\times \k(0))\stackrel{\circ^{m,0}_{p-j+1}\times {\rm Id}}{\to}
\k(l)\times\k(m-1)\stackrel{\circ^{l-1,m}_j}{\to}\k(l+m-2)$  if
$1\leq p-j+1\leq m$; otherwise define $\circ^{n,0}_p$ as
$(\k(l)\times \k(0))\times \k(m)\stackrel{\circ^{l,0}_\sigma\times {\rm Id}}{\to}
\k(l-1)\times \k(m)\stackrel{\circ^{l-1,m}_{j'}}{\to}\k(l+m-2)$,
where
$\sigma=p$, $j'=j-1$ if $p<j$; $\sigma=p-m+1$, $j'=j$ if $p>j+m-1$. 
The induction assumption implies that thus constructed maps agree on the
intersections of the faces and define a map $\partial\k(n)\to \k(n-1)$.

ii) continue $\circ^{n,0}_j$ from $\partial\k(n)$ to $\k(n)$ by sending
the center $O_n\in \k(n)$ to the center $O_{n-1}\in\k(n-1)$ and by linear 
continuation on each element of the decomposition of $\k(n)$.
Let us check that $\circ^{n,0}_j$ satisfy the induction assumption.
We need to check the operadic identities. We have two cases.

a) \begin{equation}
   \begin{matrix}
(\k(i)\times\k(0))\times \k(j) & \longrightarrow &\k(i-1)\times \k(j)& \to & \k(i+j-2)\\
\downarrow                    &                  &                  & \nearrow    &     \\
(\k(i)\times\k(j))\times \k(0)&\longrightarrow  &\k(i+j-1)\times\k(0)&&
\end{matrix}
\end{equation}

b)\begin{equation}
   \begin{matrix}
\k(i)\times(\k(j)\times \k(0)) & \to &\k(i)\times \k(j-1)& \to & \k(i+j-2)\\
\downarrow                    &     &                  & \nearrow    &     \\
(\k(i)\times\k(j))\times \k(0)&\to  &\k(i+j-1)\times\k(0)&    &
\end{matrix}
\end{equation}
We need to check that these diagrams are commutative. If $i,j<n$, then the 
commutativity follows immediately. The only cases that are left are 

1) $i=n$, $j=0$. We have two maps $\k(n)\times \k(0)\times \k(0)\to \k(n-2)$.
The induction assumption implies that they coincide on $\partial \k(n)$ and
that they map $O_n$ to $O_{n-2}$. Since both maps are piecewise linear, they coincide.
2) $i=n$, $j=1$. This is obvious.
\subsubsection{Coalgebra structure on $C_{\bullet}(\k)$} Let $C_{\bullet}(\k(n))$
be the chain complex of $\k(n)$ with respect to the polyhedral decomposition.
Since the operadic maps are piecewise linear, these complexes form a dg-operad.We are going to endow this operad with a structure of an operad of coalgebras.
To do it we need some preparation. In  sections
\ref{limits}-\ref{end} all differentials have degree -1.

\subsubsection{}\label{limits} Let $dgcoalg_1$ be the category of dg coalgebras with counit.
Then all small direct limits exist in it and commute with the forgetful functor 
$dgcoalg_1\to complexes$
\subsubsection{Cone over coalgebra} Define a coalgebra $C_{\bullet}(I)$
such that $C_0(I)$ is spanned by 2 elements $a$ and $b$; $C_1(I)$ is 
spanned by an element $c$, and
$da=db=0$; $dc=a-b$; $\Delta a=a\otimes a$; $\Delta b=b\otimes b$;
$\Delta c=c\otimes a+b\otimes c$; $\varepsilon(a)=\varepsilon(b)=1$; $\varepsilon(c)=0$.For a dg coalgebra with counit $A$ set ${\rm Cyl}\; A=
 C_\bullet(I)\otimes A$
and ${\rm Con}\; A$ to be the limit of the diagram
\begin{equation}
\begin{matrix}
{\rm Cyl}\;A &        & \\
\uparrow\lefteqn{ x\mapsto b\otimes x} &  &\\
A           &\stackrel{\varepsilon}{\longrightarrow}& k 
\end{matrix}
\end{equation}
As a vector space ${\rm Con}\; A\cong A\oplus A[-1] \oplus k$. Denote 
by $T:A\to A[-1]$ the canonical map. Denote by $e$ the element $1\in k\subset
{\rm Con}\; A$. For $u\in A$ we have:
\begin{eqnarray*}
\Delta_{\cc A}u &=& \Delta_Au\\
\Delta_{\cc A}Tu&=& (T\otimes{\rm Id} )\Delta_Au+e\otimes Tu;\\
\Delta_{\cc A}e &=& e\otimes e;\\
     d_{\cc A}u &=& du;\\
    d_{\cc A}Tu &=& u-Tdu-\varepsilon(u)e;\\
     d_{\cc A}e &=& 0.
\end{eqnarray*}
For a morphism of coalgebras with counit $\phi:A\to \cc B$ define a map
$C\phi:\cc A\to \cc B$
as follows. Let $\phi:A\to B\oplus B[-1]\oplus k\cong {\rm Con}\;B$ be defined
by its components $f:A\to B$; $Tg:A\to B\to B[-1]$; $h:A\to k$.
Then $C\phi:A\oplus A[-1]\oplus k\to B\oplus B[-1]\oplus k$
is defined by the
matrix
\begin{equation}
\left(\begin{matrix}
f        & 0 & 0\\
Tg &\tilde f & 0\\
h       &  0      & {\rm Id}
\end{matrix}\right)
\end{equation},
where $\tilde f:A[-1]\to B[-1]$ is the induced map.
Note that if $A$ and $B$ are the chain complexes
of polyhedral complexes, then the map $C\phi$ corresponds
to the piecewise-linear continuation of $\phi$ such that
the vertex of $\cc A$ goes to the vertex of $\cc B$
\begin{Proposition} $C\phi$ is a homomorphism of dg  coalgebras with counit.
\end{Proposition}
\begin{pf} Preserving of counit is clear. To prove that $C\phi$ is a 
homomorphism, let us write explicitly the condition that $\phi$
is a homomorphism. We have
\begin{eqnarray*}
\Delta_{\cc B}\phi(u) &=& \Delta_{\cc B}(f(u)+Tg(u)+h(u))\\
                      &=&
\Delta_Bf(u)+(T\otimes {\rm Id})\Delta_Bg(u)+e\otimes Tg(u)\\
                      &+& h(u)e\otimes e.
\end{eqnarray*}

$$
\phi\otimes\phi(\Delta_Au)=(f+T\circ g+h)\otimes (f+T\circ g+h)\Delta_Au 
$$

The equality $\Delta_{\cc B}\phi(u)=\phi\otimes\phi(\Delta_Au)$
yields:
$$
f\otimes f\Delta_Au=\Delta_Bf(u);
g\otimes f\Delta_Au=\Delta_Bg(u);
$$
$$
h\otimes f\Delta_Au=0;
f\otimes  g\Delta_Au=0;
$$
$$
g\otimes g\Delta_Au=0;
h\otimes g \Delta_Au=e\otimes g(u);
$$
$$
f\otimes h\Delta_Au=0; g\otimes h\Delta_Au=0;
$$
$$
h\otimes h\Delta_A u =h(u)e\otimes e.
$$
1)For $u\in A\subset \cc A$ we have $C\phi\otimes C\phi(\Delta_{\cc A}u)=
\phi\otimes\phi (\Delta_Au)=\Delta_{\cc B}C\phi(u)$  

2) For $Tu\in A[-1]\subset \cc A$
\begin{multline}
 C\phi\otimes C\phi(\Delta_{\cc A}Tu)=
(C\phi\otimes C\phi)(T\otimes {\rm Id}
\Delta_Au+e\otimes Tu)\\
=(\tilde f\otimes (f+Tg+h))(T\otimes {\rm Id}\;)\Delta_Au+(1\otimes \tilde f)(e\otimes Tu)\\
=(Tf\otimes(f+Tg+h))
 \Delta_Au+(1\otimes Tf)(e\otimes u)\\
=(T\otimes {\rm Id})\Delta_Bf(u)+(1\otimes Tf)(e\otimes u).
\end{multline}

$$
\Delta_{\cc B}C\phi(Tu)=\Delta_{\cc B}(Tfu)=(T\otimes 1)\Delta_Bf(u)+
e\otimes Tf(u).
$$
3) For 
 $e\in k\subset \cc A$, $\Delta_{\cc B} C\phi(e)=\Delta_{\cc B}e=e\otimes e$;

$C\phi\otimes C\phi\Delta_{\cc A}e=C\phi\otimes C\phi e\otimes e=e\otimes e$.

\noindent 4) The preservation of the differential and the counit
can be checked directly.
\end{pf}
\subsubsection{Category of trees} Let $\t$ be the category of planar trees with n inputs and at least one internal vertice. The morphisms are the contractions of internal edges. Let $F:\t\to Vect$ be a functor such that for $t\in \t$,
\begin{equation}\label{prd}
F(t)=\bigotimes\limits_{i\in vertices(t)} C_\bullet(K(v(i))),
\end{equation} where $v(i)$ is
the multiplicity of $i$ and for a contraction $c:t_1\to t_2$, $F(c)$
is defined by the corresponding operadic maps.
Here  the tensor product \ref{prd}
should be understood as a space of coinvariants 
$$
\big(\bigoplus\limits_{ \chi}
\bigotimes\limits_i C_\bullet(K(v(\chi(i)))\big)_{S_n},
$$
where $\chi$ is any isomorphism 
$\{1,\ldots, |vertices(t)|\}\to vertices(t)$.
see \cite{GK}.
 Then by the construction of associahedra, $C_{\bullet}(\partial K_n)=\lim\limits_{\stackrel{\longrightarrow}{\t}}F$.

\subsubsection{Construction of coalgebra structure}\label{end}
1) $n=0,1,2$.  $C_{\bullet}(\k(n))\cong k$ and we endow them with the standard 
coalgebra structure.

2) Suppose that we have found a structure of a dg coalgebra with counit
on $C_\bullet(\k(i))$, $i<n$ such that all operadic maps
$\circ^{i,j}_k:C_\bullet(\k(i))\otimes C_\bullet(\k(j))\to C_\bullet(\k(i+j-1))$,
$0\leq i,j,i+j-1\leq n-1$ are dg coalgebra morphisms and the counit maps 
$C_\bullet (\k(i))\to k$ are the augmentation maps $\varepsilon|_{C_i(\k(m))}=0$;
$e(p)=1$ for any 0-dimensional vertex of the decomposition.
The functor $F$ can be considered as 
a functor $\t\to dgcoalg_1$. According to the section \ref{limits}, 
$\lim\limits_{\stackrel{\longrightarrow}{\t}} F$ is canonically a coalgebra with counit and 
we define a coalgebra structure on $C_\bullet(\k(n))$ as the one on
$\cc \lim\limits_{\stackrel{\longrightarrow}{\t}} F\cong C_\bullet(\k(n))$. Let us check the
induction assumption. All the maps $\circ^{i,j}_k$, $1\leq i,j\leq n$, 
$i+j=n+1$ are automatically homomorphisms. The maps $\circ^{n,o}_j$ are produced from their restrictions 
${\rm res}_{\circ^{n,0}_j}:C_\bullet(\partial\k(n))\to C_\bullet(\k(n-1))\cong \cc(C_\bullet(\partial\k(n-1)))$ in the following way:
$ \circ^{n,o}_j=C{\rm res}_{\circ^{n,0}_j}$. Therefore, since 
${\rm res}_{\circ^{n,0}_j}$ are coalgebra morphisms by the induction assumption, so are
$\circ^{n,o}_j$.
We only need to check that the counit map is the augmentation map, which is obvious.  
  
\subsubsection{} We define an operad of dg coalgebras with counit $\a$ 
by setting $\a(n)^i=C_{-i}(\k(n))$.Obviously, the differential in this
operad has grading 1. Denote by $\v(n)$  the subspace in $\a(n)$ spanned
by the images of the polyhedra of the decomposition of $\k(n)$ which
do not lie on the boundary of $\k(n)$. 
\begin{Proposition}\label{aa} 
\begin{enumerate}
\item[1] $\a$ is free as an operad of graded vector spaces and is freely 
generated
by the subspaces $\v(n)$ ;
\item[2]  the counit map $\varepsilon:\a\to \as$ is a quasiisomorphism 
of operads of dg coalgebras with counit.
\end{enumerate}
\end{Proposition}
\subsection{An operad $\F$ and its cohomology}\label{cohom1}
We are going to define $\F$, a map $\F\to \bi$, and to prove that
the through map $\F\to \bi\to \e$ is a quasiisomorphism. 
We will construct the map $s:\hlie\{1\}\to \F$ in section \ref{last} which
will complete the proof of theorem \ref{third}. 
\subsubsection{}
We define $\F=\O(\a)\{1\}$ and $\G=\O(\a)$.  The following proposition
is an easy corollary of Propositions \ref{liberte} and \ref{aa}.
\begin{Proposition} 
An operad $\G$ is freely generated by operations
$\phi(v)^1_{k_1\ldots k_n}$, $k_1,\ldots,k_1\geq 1$  and $D^1_k$, $k>1$.
The map $\varepsilon:\a\to\as$ produces maps $t':\G\to\b$ and $t:\F\to\bi$.

\end{Proposition}
\subsubsection{} It is clear that $H^\bullet (\F)=H^\bullet(\G)\{1\}$.
Therefore it suffices to compute $H^\bullet (\G)$. 
We have through maps  $r:\F\to\bi\stackrel{k}{\to}\e$ and
$r': \G\to\b\to\e\{-1\}$.
\begin{Proposition}
The maps $r,r'$ are surjective on the level of cohomologies
\end{Proposition}
\begin{pf} It suffices to prove these statements for the corresponding
map $t'_*:H^\bullet(\G(2))\to H^{\bullet}(\e\{-1\})$. But $\G(2)\cong \b(2)$
and $k:\b(2)\to \e\{-1\}(2)$
is a quasiisomorphism by Theorem \ref{third} 1.
\end{pf}

Therefore, to prove that $r,r'$ are quasiisomorphisms, it suffices to give
 upper bounds for $H^\bullet (\G)$. This can be achieved by means of spectral sequences.
\subsubsection{Filtrations on $\G$ and $P(\a)$} 
For the definition of $P(\a)$ see Proposition \ref{PROP}.

Introduce the following filtration on $P(\a)$: $F^0P(\a)(n,m)=P(\a)(n,m)$; $F^1(P(\a)(n,m))=0$ if
$n\geq m$ and $F^1(P(\a)(n,m))=P(\a)(n,m)$ if $n>m$. $F^k\p(n,m)$ is defined as
a span of compositions of $\geq k$ operations among which at least $k$
belong to $F^1\p$. Note that since $\p$ is a PROP generated by an operad,
$\p(n,m)=0$, $n<m$ and $\p(n,n)=k[S_n]$. 

The filtration $F$ induces a filtration on $\G$ which is just a filtration
by the number of internal vertices of trees presenting elements of $\G$. It is easy to see that such a filtration is preserved by the differential on $\p$.
Before consideration of the spectral sequence of $\G$ associated with $F$
we will prove some Lemmas.

\subsubsection{$\phi_n(a)$ up to $F^2\p$} Let $V$ be a $\p$-algebra 
(or $\G$-algebra, which is the same); $x_1,\ldots,x_n\in T^{\geq 1}V$;
$a\in \a(n)$. Let $\phi_n^k:\a(n)\otimes TV^{\otimes n}\to TV\to V^{\otimes k}$ 
be the corestriction. Then 
$$
\phi_n^k(a, x_1,\ldots, x_n)=(\phi^1_n)^{\otimes k}
(\Delta_k a_1,\Delta_k x_1\ldots,\Delta_k x_n).
$$
Let $\Delta_k^{i_1\ldots i_k}:TV\to TV^{\otimes k}\to T^{i_1}V\otimes\ldots\otimes T^{i_k}V$ be the corestriction. Then
\begin{equation}\label{formula}
\phi_n^k(a, x_1,\ldots, x_n)=
\sum\limits_{i^1_1\ldots i^n_k}(\phi^1_n)^{\otimes k}
(\Delta_k a_1,\Delta_k^{i_1^1\ldots i_k^1} x_1\ldots,
\Delta_k^{i_1^n\ldots i_k^n} x_n).
\end{equation}

A summand in \eqnref{formula} belongs to $F^s\p$ if among the numbers 
\begin{equation}\label{no}
\sum\limits_{r=1}^n i_1^r,\ldots, \sum\limits_{r=1}^n i_k^r
\end{equation}

at least 
$s$ are greater than 1. 
Note that if some of $\sum\limits_{r=1}^n i_q^r$
is equal to 0, then the corresponding summand is equal to zero. Therefore,
to compute \eqnref{formula} up to $F^2\p$ one has to take the sum only
over such $i_1^1\ldots i_k^n$ that the numbers \eqnref{no} are all equal to 1
except, maybe, one of them, say
$\sum\limits_{r=1}^n i_s^r$. The corresponding summand in \eqnref{formula}
is equal to 
$$
\sigma=\sum\limits_{j_0\ldots j_n}\bigotimes\limits_{l=1}^{k}
\phi_n(\Delta^{j_0l}_ka,
\Delta^{j_1l}_kx_1,\ldots,\Delta^{j_nl}_k x_n),
$$
where
$$
\Delta_ka=\sum\limits_j\bigotimes\limits_{l=1}^k\Delta_k^{jl}a;
$$
$$
\Delta_k^{i_1^p\ldots i_k^p} x_p=\sum\limits_j\bigotimes\limits_{l=1}^k\Delta_k^{jl}x_p.
$$

For $l\neq s$, $\Delta^{j,l}_k x_r\in T^0V$ for all $r$
except one, say $r_l$. One may assume that $\Delta^{jl}_k{x_r}=1$
for $r\neq r_l$, $l\neq s$. We will have 
\begin{multline}
\sigma=\sum\otimes_{l=1}^{s-1}\phi_n(\Delta^{j_0l}_ka,1,1\ldots,\Delta_k^{j_{r_l}}x_{r_l}\ldots 1)\\
\otimes\phi_n(\Delta^{j_0s}_ka,\Delta^{j^1s}_kx_1,\ldots,
\Delta^{j_ns}x_n)\\
\otimes
\bigotimes\limits_{l=1}^{s-1}\phi_n(\Delta^{j_0l}_k a,1,1\ldots,\Delta_k^{j_{r_l}}x_{r_l}\ldots 1).\\
\end{multline}
But $\phi_n(a,1,1,\ldots,\xi,1,\ldots,1)=
\varepsilon(a)\xi$, where $\xi\in V$, $x\in \a(n)$ and
$\varepsilon\otimes\varepsilon\otimes\ldots\otimes 1\otimes\varepsilon\otimes\ldots\otimes\varepsilon(\Delta_na)=a$. Therefore we have
\begin{multline}
 \sigma=\sum\bigotimes\limits_{l=1}^{s-1}\Delta_k^{j_{r_l}l}x_{r_l} \otimes\phi_n(a\otimes\Delta_k^{j_1s}x_1\otimes\ldots\otimes\Delta_k^{j_ns}x_n )
\\
\otimes\bigotimes\limits_{l=s+1}^{k}\Delta_k^{j_{r_l}l}x_{r_l}
\end{multline}
Let $y_i\in T^{k_i}V$. Set $\phi'_n(a,y_1,\ldots,y_n)=\phi_n(a,y_1,\ldots,y_n)$
if $k_1+\ldots+k_n>1$ and $\phi'_n(a,y_1,\ldots,y_n)=0$ otherwise.
One sees that the sum of such sigmas is equal to
\begin{multline}\label{Gg}
\sum\limits_{j_1\ldots j_n}\pmin x_1^{j_11}\cup x_2^{j_21}\cup\ldots\cup x_n^{j_n1}
\\ 
\otimes \phi'_n( a, x_1^{j_12}, x_2^{j_22},\ldots, x_n^{j_n2})\\
\otimes x_1^{j_13}\cup x_2^{j_23}\cup\ldots\cup x_n^{j_n3}+\varepsilon(a)x_1\cup\ldots\cup x_n,
\end{multline}
where
$\Delta_3x_i=\sum\limits_{j_i}x_i^{j_i1}\otimes x_i^{j_i2}\otimes x_i^{j_i3},
$
and $\cup$ means shuffle.
Thus we have proven
\begin{Lemma}
$\phi_n(x)$ is equal to \eqnref{Gg} up to operations in $F^2\p$.
\end{Lemma}
\subsubsection{Additional filtration on $F^1\G/F^2\G$} Note that
$\w=F^1\G/F^2\G$ is generated by the images of the elements
$\phi(v)^1_{n_1\ldots n_k}$, $v\in \v(k)$ and $D^1_k$. Define a filtration
$F'$ on $\w$ as follows. $F'_m\w$ is generated by $\phi(v)_{n_1\ldots n_l}$,
$l\leq m$ and $D^1_l$, $m\geq 1$.
\begin{Lemma}\label{dva}
\begin{enumerate}
\item[1)] Filtration $F'$ is compatible with the induced differential on $\w$.
\item[2)] $d\phi(v)^1_{n_1\ldots n_l}=\phi(\delta v)^1_{n_1\ldots n_l}$
in ${\rm Gr}_{F'}\w$, where $\delta$
 is the differential on $\v(l)\cong C_{-\bullet}(\k(k),\partial\k(k))$.
\end{enumerate}
\end{Lemma}
\begin{pf} Take $v\in V(k)$ corresponding to an element of decomposition of $\k(k)$ denoted by the same letter. Since $v$
does not belong to $\partial \k(k)$, $v={\rm Con}\;w$, $w\subset \partial \k(k)$, $w=\circ_i^{p,k-p}(\lambda, \mu)$, $\lambda,\mu\in \k(<k)$.
We have $dv=\bar\delta v+w$, where $d$ is the differential on $\a$ and
$\bar\delta v$ is the oriented sum of all faces of $v$ that do not lie
on $\partial \k(k)$.
According to \eqnref{differen}, for any $\p$-algebra, we have
\begin{multline}\label{differe}
 d\phi^1_{k_1\ldots k_n}(v, x_1,\ldots, x_n)-  
\phi^1_{k_1\ldots k_n}(dv, x_1,\ldots, x_n)\\
 -\sum(-1)^{|v|+|x_1|+\ldots+|x_{l-1}|}\phi^1_{k_1\ldots k_n}(v, x_1,\ldots dx_l,\ldots, x_n)\\
=-\sum\limits_s D^r_s\phi^s_{k_1\ldots k_n}(v, x_1,\ldots, x_n)
 +\sum(-1)^{|v|+|x_1|+\ldots+|x_{l-1}|}\phi^1_{k_1\ldots s\ldots k_n}(v, x_1
,\ldots D_l^s x_l,\ldots, x_n)\\
\end{multline}

Hence, up to $F^2\p$, we have
\begin{multline}\label{differential}
d\phi^1_{k_1\ldots k_n}(v, x_1,\ldots, x_n)  
 -\sum(-1)^{|v|+|x_1|+\ldots+|x_{l-1}|}\phi^1_{k_1\ldots k_n}(v, x_1,\ldots,
 dx_l,\ldots, x_n)\\
=\phi^1_{k_1\ldots k_n}(dv, x_1,\ldots, x_n)
=\phi^1_{k_1\ldots k_n}(\bar\delta v, x_1,\ldots, x_n)\\
+\phi^1(\lambda, x_1,\ldots x_{i-1},\phi(\mu, x_i,\ldots, x_j),\ldots, x_n)
\end{multline}
The last expression is equal to
$$
\phi^1_{k_1\ldots k_n}(\bar\delta v, x_1,\ldots, x_n)
$$
modulo $F'_{n-1}\w$
\end{pf}
\subsubsection{Fundamental classes $\mu_k$.}\label{fund}
Let $\mu_k\in C_{k-2}(\k(k))$ be the fundamental class. It is clear that 
$\mu_k\in \v(k)$ and that
 \begin{equation}\label{difmu}
d\mu_k=\sum\limits_{i=2}^{k-1}\mu_i\{\mu_{k+1-i}\},
\end{equation}
where
$$
\mu_i\{\mu_j\}=\sum\limits_{k=1}^{i}(-1)^{(k-1)(j-1)}\circ^{i,j}_k(\mu_i,\mu_j)
$$
 For a map $\chi:TV^{\otimes n}\to V$ define
$T\chi:TV^{\otimes n}\to TV$ by the formula (cf. (\ref{Gg})

\begin{multline}\label{G}
T\chi(x_1,\ldots, x_n)=\sum\limits_{j_1\ldots j_n}\pmin x_1^{j_11}\cup x_2^{j_21}\cup\ldots\cup x_n^{j_n1}
\\ 
\otimes \chi(x_1^{j_12}, x_2^{j_22},\ldots, x_n^{j_n2})\\
 \otimes x_1^{j_13}\cup x_2^{j_23}\cup\ldots\cup x_n^{j_n3}.
\end{multline}
For $k>2$ denote by $\phi^1_k$ the corestriction of $\phi(\mu_k)$ on $V$.
For $k=2$ let $\phi^1(\mu_2)$ be the restriction-corestriction
of $\phi_2(\mu_2):T^{\geq 1}V\otimes T^{\geq 1}V\to TV\to V$.
Then we have 
\begin{equation}\label{shuffle}
\phi_2(\mu_2)(x,y)=x\cup y+T\phi^1(\mu_2)(x,y)+k(x,y),
\end{equation}
where $k\in F^2\p.$
Finally, let $\phi_1$ be the restriction-corestriction of the differential
$D:T^{\geq2}V\to TV\to V$.
\begin{Lemma}\label{tri}
\begin{equation}
d\phi_k(x_1,\ldots,x_k)=\sum\limits_{i=1}^k\phi_i\{T\phi_{k+1-i}\}+
\sum\limits_r(-1)^{r-1}\phi_{k-1}(x_1,\ldots x_r\cup x_{r+1},\ldots)+u,
\end{equation}
where $u\in F^3\G$ and 
\begin{multline*}
\phi_i\{T\phi_{k+1-i}\}(x_1,\ldots,x_k)\\
=\sum\limits_{l=1}^i(-1)^{(l-1)(k-i)}
(-1)^{(k+1-i)(|x_1|+\ldots+|x_{l-1}|)}\phi_i(x_1,\ldots,x_{l-1},T\phi_{k+1-i}(x_l,\ldots,x_{l+k-i}),\ldots,x_k)
\end{multline*}
\end{Lemma}
\begin{pf} This immediately follows from \eqnref{Gg}, \eqnref{difmu}, and
\eqnref{shuffle}
\end{pf}
\subsection{Spectral Sequence Associated with $ F$}\label{cohom2}
\subsubsection{$E_1^{\bullet,\bullet}$} Since $\G$ is free, the K\"unneth
formula implies that the operad of $E_1^{\bullet,\bullet}$ is a 
semi free operad generated by $H^\bullet(\w)$. Let $\z\subset \w$ be
the $S$-submodule generated by the images of $\phi(\mu_k)^1_{n_1\ldots n_k}$.
\begin{Lemma}
\begin{enumerate}
\item[1)] $\z$ is a subcomplex of $\w$.
\item[2)] ${\rm Gr}_{F'}\;\z\to {\rm Gr}_{F'}\;\w$
is a quasiisomorphism.
\item[3)] $\z \to \w$ is a quasiisomorphism
\end{enumerate}
\end{Lemma} 
\begin{pf}
1). Follows from Lemma \ref{tri} which implies that
\begin{equation}\label{diffsh}
d\phi_k(x_1,\ldots,x_k)=\sum\limits_r(-1)^{r-1}\phi_{k-1}(x_1\otimes\ldots x_r\cup x_{r+1}\otimes\ldots)+F^2\G.
\end{equation}

2) Follows from Lemma \ref{dva}.

3) Follows from 2).
\end{pf}
 
Consider the complex $\z$. It is clear that the Schur functor
$S_{\z}V\cong\oplus_{n=2}^{\infty}\z(n)\otimes_{S_n}T^nV\cong
\bigoplus\limits_{k=1}^{\infty}T^k(T^{\geq 1}V)[k-2]\cong T(T^{\geq 1}V[1])[-2]$
for any vector space $V$. The differential \ref{diffsh} coincides
with the Hochschild differential of the algebra $(T^{\geq 1}V,\cup)$.
Therefore, 
$$
S_{H^\bullet (\z)}V\cong HH_{2-\bullet}(T^{\geq 1}V)\cong
$$
$$
 S^\bullet(T^{\geq 1}V/shuffles[1])[-2]\cong P(V)[-1]\cong
S_{(\e\{-1\})^*[-1]}V.
$$
An $S$-module $\g=(\e\{-1\})^*[-1]$ constitutes the space of generators
of the canonical resolution $\he\{-1\}$

Therefore,  $E_1^{\bullet,\bullet}$
is a free operad generated by $\g$.
We have a quasiisomorphism 
\begin{equation}\label{aph}
\alpha:\z\to H^\bullet(\z)
\end{equation}
such that
the map $\alpha_*:S_\z V\to S_{H^\bullet(\z)}V$ is the obvious projection
$ T(T^{\geq 1}V[1])[-2]\to S^\bullet(T^{\geq 1}V/shuffles[1])[-2]$

 We are going to prove that the differential 
$d_1$ coincides with the differential on $\he\{-1\}$.
This implies  that $E_2^{\bullet,\bullet}\cong \e\{-1\}$ 
and that the map $r':\G\to\e\{-1\}$ is a quasiisomorphism.
\subsubsection{$ d_1$.} This differential is a quadratic differential
on a free operad $F(\g)$
generated by $\g$.
Let $\g(n)=\e\{-1\}(n)^*[-1]$.
 The space of quadratic elements of $F(\g)$
is isomorphic to 
$\bigoplus \g(i)\otimes 
\g(j)\otimes_{S_j\times S_{i-1}}k[S_{i+j-1}]$, where we assume that the second factor
is inserted into the first position of the first factor. Any such a differential defines on $\g[1]=\e\{-1\}^*$ a structure of cooperad. Let us denote the cooperad
corresponding to $d_1$ by $e'$.  We need to prove that $e'= \e\{-1\}^*$. 

Note that $\e\{-1\}^*$ is cogenerated by $\e\{-1\}^*(2)$. Therefore,
it suffices to prove that the coinsertion maps 
$o_1^{2,n*}:e'(n+1)\to e'(2)\otimes e'(n)$ are the same as in $\e\{-1\}^*$.
For this we need to compute the restriction ${\bar d}_1:\g(j+1)\to\g(2)\otimes\g(j)\otimes_{S_j}k[S_{j}]$.
Take an operation $\xi=\varphi(\mu_l)^1_{k_1\ldots k_n}$. Then (\ref{difmu})
allows one to compute ${\bar d}_1\alpha(\xi)$,
where $\alpha$ is the map (\ref{aph}).
 After removing
 irrelevant terms, we obtain
\begin{multline*}
{\bar d}_1\alpha(\xi(x_1,\ldots, x_l))=
\alpha\Big(D_2(T\phi_l(x_1,\ldots, x_l))+
\phi(\mu_2)_{1,1}(x_1,\phi_{l-1}(x_2\ldots x_l))\\
\pmin 
\phi(\mu_2)_{1,1}(\phi_{l-1}(x_1\ldots x_{l-1}),x_l\Big).
\end{multline*}
Let
$\Delta x_i=\sum\limits_{pq}\Delta^{p,q} x_i$, where 
$\Delta^{p,q} x_i\in T^pV\otimes T^qV$. Let
 $\Delta^{p,q}x_i=\sum\limits_t  x_{i1}^{pqt}\otimes x_{i2}^{pqt}$. Then
\begin{multline*}
{\bar d}_1\alpha(\xi(x_1,\ldots, x_l))=
\alpha\Big(\sum D_2(x^{1,q}_{r1}\otimes\phi^1(x_1,\ldots ,x^{1,q}_{r2},\ldots ,x_l))\\
+\sum D_2(\phi^1(x_1,\ldots, x^{q,1}_{r1},\ldots, x_l)\otimes x^{q,1}_{r2})
\\
+\phi(\mu_2)_{1,1}(x_1,\phi_{l-1}(x_2\ldots x_l))\pmin 
\phi(\mu_2)_{1,1}(\phi_{l-1}(x_1\ldots x_{l-1}),x_l)\Big).
\end{multline*}
Note that $\alpha(D_2(x\otimes y))$ is the 
cocommutator in $\e\{-1\}^*$
and $\alpha(\phi(\mu_2)_{1,1}(x,y))$ is the coproduct.
Now the coincidence of the cooperadic structures on $e'$
and $\e\{-1\}^*$ is obvious.

\subsection{Maps $s:\hlie\{1\}\to \F$ and $s':\hlie\to\G$}\label{last} We construct a map $s
$ satisfying the Theorem \ref{third}. We have a map of dg operads
$\hass \to \a$ such that $m_k\in \hass\mapsto \mu_k\in \a$,
where $\mu_k$ are the fundamental classes, see (\ref{fund}).  
Also we have a map $\hlie\to\hass:m_k({x_1\ldots x_k})\to \sum\limits_\sigma(-1)^\sigma\mu_n( x_{\sigma_1}\ldots x_{\sigma_n})$.
Therefore, any $\O(\a)=\G$-structure on $V$ implies a structure of a $\hlie$-
algebra on $TV$. The operations are given by the antisymmetrization
of $\phi(\mu_k)(x_1\ldots x_k)$.
\begin{Lemma}The operations $\phi^r(\mu_k)(x_1\ldots x_k)$, $k\geq 3$
and $\phi^r(\mu_2)(x_1\otimes x_2)\mplus\phi^r(\mu_2)(x_2\otimes x_1)$
are equal to zero whenever all $x_i$ are in $T^1V$ and $r\geq 2$.
\end{Lemma}\begin{pf}1) $k\geq 3$. We have 
\begin{multline}
\phi^r(\mu_k)(x_1\ldots x_k)\\
=(\phi^1)^{\otimes r}(\Delta_r\mu_k,\Delta_rx_1,\ldots,\Delta_rx_k)\\
=\sum\bigotimes\limits_{i=1}^r\phi^1(\Delta_r^{j_0i}\mu_k,
\Delta_r^{j_1i}x_1,\ldots,\Delta_r^{j_ki}x_k).
\end{multline}  This  is a sum of products. Take one of these products.
Let the $i$-th multiple in it have $k-s_i$ arguments to be equal to 1
and $s_i$ arguments from $V$. We have $\sum s_i=k$.  This  product can be rewritten as a product of operations $\phi(a_i)$, where $a_i\in \k(s_i)$.
Since ${\rm dim}\;\k(i)=i-2$ if $i\geq 2$ and $\dim \k(1)=0$, the total 
grading
of such an operation is at least $-(s_1-1+\ldots+s_r-1)=-(k-r)$. Since the degree
of $\mu_{k}$ is $-(k-2)$, the product vanishes if $r\geq 3$. If 
$r=2$, then the equality of the degrees is only possible when $s_1=s_2=1$,
that is when $k=2$.

2)$k=2.$ The above argument allows us to assume $r=2$.
We have $\phi^2(\mu_2\otimes x\otimes y)\mplus\phi^2(\mu_2\otimes x\otimes y)=
\phi^1\otimes \phi^1((\mu_2\otimes\mu_2),([x,y]),(1\otimes 1) )-
\phi^1\otimes \phi^1((\mu_2\otimes\mu_2),(1\otimes 1),([x,y]))=0.$
\end{pf}
Hence, the $\hlie$-algebra on $TV$ can be restricted to $V$. Therefore,
we have a map $s':\hlie\to\G$. To compute the through map 
$\hlie\{1\}\to\F\to\bi$, let us notice that the image of $\mu_k$
under the counit map $\a\to\as$ is zero for $k>2$. Therefore all higher
products will go to zero, and $\phi^1(\mu_2, x, y)\mplus  
 \phi^1(\mu_2, y, x)$ will go to $m_{1,1}(x,y)\mplus m_{1,1}(y,x)$,
which proves that the diagram in Theorem \ref{third} commutes.
     
\end{document}